\documentclass[a4paper, 12pt]{article}

\usepackage{amsmath, amssymb, latexsym, amsthm, amsfonts, setspace, bm}
\usepackage[dvips]{color, graphicx}
\usepackage{float}

\theoremstyle{definition}
\newtheorem{definition}{Definition}[section]
\newtheorem{thm}[definition]{Theorem}

\newtheorem{cor}[definition]{Corollary}
\newtheorem{lem}[definition]{Lemma}

\newtheorem{rem}[definition]{Remark}

\newtheorem*{question}{Question}
\newtheorem*{note}{Note} 

\theoremstyle{definition}
\numberwithin{equation}{section}

\begin{document}
\title{A sufficient condition for a graph with boxicity at most its chromatic number}
\author{\footnote{Academic Assembly Institute of Science and Engineering, Shimane University, Matsue, 
Shimane 690-8504, Japan. \endgraf 
\hspace{0.27cm}{\it E-mail address}:\,kamibeppu@riko.shimane-u.ac.jp}\,\,Akira Kamibeppu}
\date{}
\maketitle
\setstretch{1.0}

\begin{abstract}
A {\it box} in Euclidean $k$-space is the Cartesian product of $k$ closed intervals on the real line.\ 
The {\it boxicity} of a graph $G$, denoted by $\text{box}(G)$, is the minimum nonnegative integer $k$ 
such that $G$ can be isomorphic to the intersection graph of a family of boxes in Euclidean $k$-space.\ 
In this paper, we present a sufficient condition for a graph $G$ under which $\text{box}(G)\leq \chi (G)$ 
holds, where $\chi (G)$ denotes the chromatic number of $G$.\ 
Bhowmick and Chandran \cite{BC10} proved that $\text{box}(G)\leq \chi (G)$ holds for a graph $G$
with no asteroidal triples.\ We prove that $\text{box}(G)\leq \chi (G)$ holds for a graph $G$ in a special 
family of circulant graphs with an asteroidal triple.\\

\vspace{0.1cm}
\end{abstract}
\maketitle

\section{Introduction and Preliminaries}
A {\it box} in Euclidean $k$-space is the Cartesian product of $k$ closed intervals on the real line.\
The {\it boxicity} of a graph $G$, denoted by $\text{box}(G)$, 
is the minimum nonnegative integer $k$ such that $G$ can be isomorphic to the intersection graph of 
a family of boxes in Euclidean $k$-space.\ The concept of boxicity of graphs was introduced 
by Roberts \cite{Ro69}.\ It has applications in some research fields, for example, a problem of niche overlap 
in ecology (see \cite{Ro76} for detail).\ 
It is known that boxicity of graphs have relationships with some graph invariants.\ 
Chandran, Francis and Sivadasan \cite{CFS08} proved that $\text{box}(G)\leq 2\Delta _G^2$ 
holds for a graph $G$, where $\Delta _G$ is the maximum degree of $G$.\ In fact, they also showed 
in the proof that $\text{box}(G)\leq 2\chi (G^2)$ holds for a graph $G$, where $G^2$ is the graph 
obtained from $G$ by adding edges whose endvertices have common neighbors in $G$ and 
the symbol $\chi (G^2)$ means the chromatic number of $G^2$.\  Moreover they conjectured 
that $\text{box}(G)$ is $O(\Delta _G)$.\ 
Esperet \cite{E09} improved the previous upper bound for boxicity, that is, proved 
that $\text{box}(G)\leq \Delta _G^2+2$ holds for a graph $G$.\ 
Adiga, Bhowmick and Chandran \cite{ABC11} disproved the above conjecture.\ In fact, they 
proved that there exist graphs with boxicity $\Omega (\Delta _G \log \Delta _G)$.\ Before these 
results appear, Chandran and Sivadasan \cite{CS07} presented chordal graphs, circular arc graphs, 
asteroidal triple free graphs, co-comparability graphs and permutation graphs as examples whose boxicity 
are bounded above by a linear function of their maximum degree.\ 
These results are based on their main result which states $\text{box}(G)\leq \text{tw}(G)+2$ for a graph $G$, 
where $\text{tw}(G)$ is the treewidth of $G$.\

So far some researchers found relationships between boxicity and chromatic number.\ 
Chandran, Das and Shah \cite{CDS09} proved that if the boxicity of a graph $G$ with $n$ vertices 
is equal to $n/2-s$ for $s\geq 0$, then $\chi (G)\geq n/(2s+2)$ holds.\ This result implies that if the 
boxicity of a graph is close to $n/2$, then its chromatic number is large.\ 
Esperet \cite{E13} proved that $\text{box}(G)\leq \chi _a(G)(\chi _a(G)-1)$ holds for a graph $G$ with 
$\chi _a(G)\geq 2$, where $\chi _a(G)$ is the acyclic chromatic number of $G$.\ 
Bhowmick and Chandran \cite{BC10} proved that $\text{box}(G)\leq \chi (G)$ holds 
for a graph $G$ with no asteroidal triples.\ We remark that $\text{box}(G)\leq \chi (G)$ does not hold in general.\ 
We consider the graph $H$ obtained from a balanced complete bipartite graph with at least 
10 vertices by removing a perfect matching.\ Then $\text{box}(H)>2=\chi (H)$ holds 
(see \cite{CDS09} for detail).\ 
Recently, Chandran Chandran, Mathew and Rajendraprasad \cite{CMR16} remarked that almost all graphs 
have boxicity more than their chromatic number, which is based on the probabilistic method, 
but the family of graphs with boxicity at most their chromatic number is not narrow.\ 
In this paper, we present a sufficient condition for a graph $G$ under which $\text{box}(G)\leq \chi (G)$ 
holds.\ Moreover we show that $\text{box}(G)\leq \chi (G)$ holds for a graph $G$ in a special family of 
circulant graphs with an asteroidal triple.\ Other results appear in \cite{CIMR15, CMS11}.\ 

All graphs are finite, simple and undirected in this paper.\ We use $V(G)$ for the vertex set of a graph 
$G$ and $E(G)$ for the edge set of the graph $G$.\ An edge of a graph with endvertices $u$ and $v$ is 
denoted by $uv$.\ Then $u$ is called a neighbor of $v$ (or $v$ is called a neighbor of $u$).\ 
The set of all neighbors of a vertex $v$ in $G$ is denoted by $N_G(v)$, or briefly by $N(v)$.\ 
A map $c:V(G)\rightarrow \{1,\ldots , k\}$ is called a {\it $k$-coloring} of a graph $G$ 
if $c(u)\ne c(v)$ holds whenever $uv$ is an edge of $G$.\ 
The {\it chromatic number} of $G$, denoted by $\chi (G)$, is the minimum positive integer $k$ such that 
the graph $G$ has a $k$-coloring.\ A $k$-coloring of $G$ gives a partition of $V(G)$ into $k$ independent 
sets, called {\it color classes}, so we use the notation $\{V_i\}_{i=1}^k$ for a $k$-coloring of $G$ with 
color classes $V_1, \ldots , V_k$.\

The following is a useful result to calculate boxicity of graphs.\ 
\begin{thm}[\cite{Ro69}]Let $G$ be a graph.\ Then $\text{box}(G)\leq k$ holds if and only if 
there exist $k$ interval graphs $H_1, \ldots , H_k$ with vertex set $V(G)$ such that 
$E(G)=E(H_1)\cap \cdots \cap E(H_k)$ holds.\ 
\end{thm}

\section{Split Interval Graphs}
For a $k$-coloring $\{V_i\}_{i=1}^k$ of a graph $G$, we will consider a supergraph $H_i$ of $G$ 
which accompanies a color class $V_i$ so that $E(G)=E(H_1)\cap \cdots \cap E(H_k)$ holds.\ 
Nonadjacent vertices of $H_i$ might have at least two common neighbors in $V_i$, which causes induced 
cycles of $H_i$, so making the set $V(H_i)\setminus V_i$ complete is a natural way to avoid their cycles.\ 
The resulting graph is called a split graph (see Theorem 2.1 below).\ 

A graph is {\it chordal} if the graph contains no induced cycles other than triangles.\  
A graph is called a {\it split} graph if the graph and its complement are chordal.\ The following 
is a characterization of split graphs.\ 
\begin{thm}[\cite{FH77}]
A graph $G$ is a split graph if and only if $V(G)$ can be partitioned into an independent set and a clique.\ 
\end{thm} 
In this paper, $P_{uv}$ denotes a path between vertices $u$ and $v$ of a graph $G$, if exists.\  
A triple of vertices $u$, $v$ and $w$ of $G$ is said to be {\it asteroidal} if there exist paths $P_{uv}$, $P_{vw}$ 
and $P_{wu}$ in $G$ such that $N_G(w)\cap V(P_{uv})=\emptyset $, $N_G(u)\cap V(P_{vw})=\emptyset $ and 
$N_G(v)\cap V(P_{wu})=\emptyset $ hold.\ 
We introduce a characterization of interval graphs.\ 
\begin{thm}[\cite{LB62}]
A graph $G$ is an interval graph if and only if $G$ is chordal and has no asteroidal triples.\ 
\end{thm}

\begin{lem}
Let $G$ be a split graph with a partition of $V(G)$ into an independent set $S$ and a clique $K$.\ 
If for any triple of vertices $u$, $v$ and $w$ in $S$ two vertices $x$ and $y$ in $\{u,v,w\}$ satisfy 
$N_G(x) \supseteq N_G(y)$, the graph $G$ is an interval graph.\ 
\begin{proof}
It is sufficient by Theorem 2.2 to show that any triple of vertices of the split graph $G$ is not 
asteroidal since the family of split graphs is a special class of chordal graphs.\ 
For any triple of vertices $u$, $v$ and $w$ of $G$, if one of them, say $u$, is in the clique $K$, then 
$N_G(u)$ contains vertices of any path $P_{vw}$.\ 
So in this case the triple of vertices is not asteroidal.\ Hence we may assume that the triple of vertices 
$u$, $v$ and $w$ of $G$ are in the independent set $S$.\ By our assumption, if 
$N_G(u)\supseteq N_G(v)$ holds, $N_G(u)$ contains a vertex in any path $P_{vw}$.\ Hence our 
assumption tells us that the triple is also not asteroidal.\ 
\end{proof}
\end{lem}

\section{Our Results}
The following is our main result in this paper.\ 
\begin{thm}
If a graph $G$ has a $\chi (G)$-coloring $\{V_i\}_{i=1}^{\chi (G)}$, where 
$V_i=\{v_{i,1}, v_{i,2}, \ldots , v_{i, n(i)}\}$, with the property that, 
for each vertex $v_{i,j}$ in $V(G)$, 
there exists a subset $X_{i,j}$ of 
$V(G)\setminus V_i$ containing $N_G(v_{i,j})$ such that 
\begin{itemize}
\item[(i)] $X_{i,1}\supseteq \cdots \supseteq X_{i,k(i)}$ and 
$X_{i,k(i)+1}\subseteq \cdots \subseteq X_{i,n(i)}$ hold for each $i$, 
where $k(i)\in \{1,2, \ldots , n(i)\}$, and 
\item[(ii)]either $v_{i,s}\not\in X_{j,t}$ or $v_{j,t}\not\in X_{i,s}$ holds for any pair of nonadjacent 
vertices $v_{i,s}$ and $v_{j,t}$ in $G$, 
\end{itemize}
then the inequality $\text{box}(G)\leq \chi (G)$.\
\begin{note}
In (i), `$k(i)=n(i)$' means that only the sequence $X_{i,1}\supseteq \cdots \supseteq X_{i,n(i)}$ is required.\ 
It also has a similar meaning in Corollary 3.3 below.\  
\end{note}
\begin{proof}
Let $\{V_i\}_{i=1}^{\chi (G)}$ be a $\chi (G)$-coloring of the graph $G$ with the property.\ 
We define the supergraph $H_i$ of $G$ for each color $i\in \{1,2,\ldots , \chi (G)\}$ as follows: 
\begin{equation*}
V(H_i)=V(G), \hspace{0.2cm} E(H_i)=E(G)\cup \bigcup_{k=1}^{n(i)}\{ v_{i,k}x\,|\, x\in X_{i,k}\}
\cup \{xy\,|\,x,y\in V(G)\setminus V_i \}.
\end{equation*}
By construction the graph $H_i$ is a split graph with the partition of $V(H_i)$ into the independent set 
$V_i$ and the clique $V(G)\setminus V_i$ for each $i\in \{1,2,\ldots , \chi (G)\}$.\ For any triple of vertices 
$v_{i,p}$, $v_{i,q}$ and $v_{i,r}$ in $V_i$, we can find two vertices $v_{i,s}$ and $v_{i,t}$ in 
$\{v_{i,p}, v_{i,q}, v_{i,r}\}$ such that $X_{i,s}\supseteq X_{i,t}$ by condition (i).\ Hence we see that 
$H_i$ is an interval graph by Lemma 2.3 since $X_{i,j}$ is the same as $N_{H_i}(v_{i,j})$.\

We note that each edge of $G$ is an edge of $H_i$ for each $i\in \{1,2,\ldots , \chi (G)\}$ by definition.\ 
Let $u$ and $v$ be nonadjacent vertices in $G$.\ We write $u=v_{i,k}$ and $v=v_{j,l}$ under the 
coloring $\{V_i\}_{i=1}^{\chi (G)}$.\  If $i=j$, then $u$ and $v$ are in $V_i$, and hence they are 
not adjacent in $H_i$.\ In what follows, we assume $i\ne j$.\ By condition (ii), $u$ and $v$ are 
nonadjacent in $H_i$ if $v_{j,l}\not\in X_{i,k}$ holds, otherwise $u$ and $v$ are nonadjacent in $H_j$ since 
$v_{i,k}\not\in X_{j,l}$ holds.\ Hence we have $E(G)=E(H_1)\cap \cdots \cap E(H_{\chi (G)})$.\ 
This completes the proof of our theorem by Theorem 1.1.\ 
\end{proof}
\end{thm}
Every optimal vertex coloring $c$ of a complete multipartite graph satisfies the condition of Theorem 3.1.\ 
For each vertex $v$ of the graph, take a set $X_{c(v), v}$ as the set of all neighbors of $v$.\ 
Paths and cycles also have an optimal coloring with the condition of Theorem 3.1.\ 
Needless to say, we see that the boxicity of these graphs are at most their chromatic number respectively 
without applying Theorem 3.1.\ 
\begin{rem}
The family of graphs with vertex colorings that satisfy the condition of Theorem 3.1 is not contained 
by the family of graphs with no asteroidal triples.\ Consider the graph in Figure 1.\ It is easy to check that 
the triple of vertices with degree 1 in the graph forms an asteroidal triple.\ Under the vertex 
coloring $c$ of the graph in Figure 1, we define $X_{c(v), v}$ as the set of all neighbors of $v$.\ 
\begin{figure}[!ht]
\centering
\includegraphics[scale=0.7,clip]{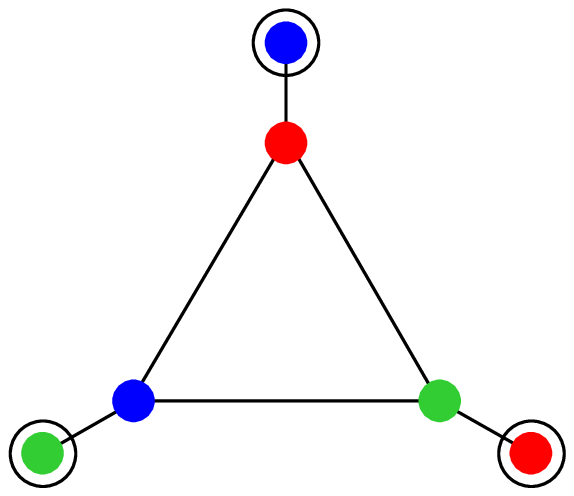}
\renewcommand{\baselinestretch}{1}
\caption[]{A graph with an asteroidal triple.}
\end{figure}
\end{rem}
The condition which the colorings of complete multipartite graphs or the graph in Figure 1 satisfy 
is stronger than the condition which Theorem 3.1 claims. 
\begin{cor}
If a graph $G$ has a $\chi (G)$-coloring $\{V_i\}_{i=1}^{\chi (G)}$ with the property that, 
for each color $i\in \{1,2,\ldots , \chi (G)\}$, 
\begin{center}
$N_G(v_{i,1})\supseteq \cdots \supseteq N_G(v_{i,k(i)})$ and 
$N_G(v_{i,k(i)+1})\subseteq \cdots \subseteq N_G(v_{i,n(i)})$ 
\end{center}
hold, where $V_i=\{v_{i,1}, v_{i,2}, \ldots , v_{i, n(i)}\}$ and $k(i)\in \{1,2, \ldots , n(i)\}$, the inequality 
$\text{box}(G)\leq \chi (G)$ holds.\ 
\begin{proof}
Define $X_{i,k}$ as $N_G(v_{i,k})$ for each vertex $v_{i,k} \in V(G)$.\ 
For any pair of nonadjacent vertices $v_{i,k}$ and $v_{j,l}$, we note that $v_{i,k}$ is not in $N_G(v_{j,l})$ 
(and $v_{j,l}$ is also not in $N_G(v_{j,l})$).\ The corollary follows from Theorem 3.1.\ 
\end{proof}
\end{cor}

\section {Boxicity of Circulant Graphs}
In this section we present a family of graphs with an asteroidal triple and show that $\text{box}(G)\leq \chi (G)$
for a graph $G$ in the family.\ 
We define the graph $G_{a,b}$, where $a\geq 2b$, as follows: $V(G_{a,b})=\{0,1, \ldots , a-1\}$ and 
$uv\in E(G_{a,b})$ if and only if $u\in \{v+b, v+b+1, \ldots , v+a-b\}$ with addition modulo $a$.\ 
For example, $G_{a,1}$ and  $G_{a,2}$ are isomorphic to a complete graph and the complement of a cycle 
respectively.\ The graph $G_{a,b}$ is a member of the family of circulant graphs.\ It is known that 
$\alpha (G_{a,b})=b$ and $\chi (G_{a,b})\geq a/b$ hold (see \cite{SU}, page 44 for detail), where 
$\alpha (G_{a,b})$ denotes the independence number of $G_{a,b}$.\ 
\begin{thm}
For a circulant graph $G_{nb, b}$ with $n\geq 2$ and $b\geq 1$, the inequality 
$\text{box}(G_{nb,b}) \leq \chi (G_{nb,b})$ holds.\ 
\begin{proof}
We arrange all vertices $0,1, \ldots, nb-1$ of $G_{nb,b}$ clockwise in this order.\ We see that every set of 
consecutive $b$ vertices forms a maximum independent set of $G_{nb,b}$.\ We also note that 
$n=nb/b\leq \chi (G_{nb,b})\leq n$, that is, $\chi (G_{nb,b})=n$ holds since $V(G_{nb,b})$ 
can be partitioned into $n$ maximum independent sets with $b$ vertices.\
 
We write $v_{i,j}$ for the vertex $(i-1)b+j-1$ in $V(G_{nb,b})$, 
where $i\in \{1,2,\ldots, n\}$ and $j\in \{1,\ldots ,b\}$.\ 
We note that $v_{0,j}$ and $v_{n+1,j}$ are identified with $v_{n,j}$ and $v_{1,j}$ in $G_{nb,b}$ respectively.\ 
Let $V_i$ be the set $\{v_{i,1}, \ldots , v_{i,b}\}$.\ 
We see that $\{V_i\}_{i=1}^n$ becomes an $n$-coloring of $G_{nb,b}$.\ 
We also see that 
\begin{equation*}
N(v_{i,j})=\{v_{i-1,1},\ldots , v_{i-1,j}\}\cup \{v_{i+1,j}, \ldots , v_{i+1,b}\}\cup U_i 
\end{equation*}
holds where $U_i=V(G_{nb,b})\setminus (V_{i-1}\cup V_i\cup V_{i+1})$.\ We consider cases 
whether $b$ is even or not.\ In either case our statement follows from Theorem 3.1. \vspace{0.2cm}\\  
\noindent {\bf Case (I)}: We assume that $b$ is even.\vspace{0.2cm}\\
We define the subset $X_{i,j}$ of $V(G_{nb,b})\setminus V_i$ for the vertex $v_{i,j}$ of $G_{nb,b}$ as follows:
\begin{equation*}
X_{i,j}=
\begin{cases}
N(v_{i,j})\cup \{v_{i-1,j},\ldots , v_{i-1, b/2}\}
& \text{if $1\leq j\leq b/2$},\\
N(v_{i,j})\cup \{v_{i+1,b/2+1},\ldots , v_{i+1,j}\}
& \text{if $b/2<j\leq b$}.
\end{cases}
\end{equation*}
We can check that  
\begin{center} 
$X_{i,1}\supseteq \cdots \supseteq X_{i,b/2}$ and 
$X_{i,b/2+1}\subseteq \cdots \subseteq X_{i,b}$
\end{center}
hold for each color $i$.\ 

We remark that the vertex $v_{i,j}$ of $G_{nb,b}$ is not adjacent to every vertex in 
\begin{center}
$\{v_{i-1,j+1}, \ldots , v_{i-1,b}\}\cup V_i \cup \{v_{i+1,1},\ldots , v_{i+1,j-1}\}$
\end{center}
for $2\leq j\leq b-1$.\ The vertex $v_{i,1}$ is not adjacent to every vertex in 
$(V_{i-1}\setminus \{v_{i-1, 1}\})\cup V_i $ and $v_{i,b}$ is not adjacent to every vertex in 
$V_i \cup (V_{i+1}\setminus \{v_{i+1,b}\})$.\ \vspace{0.2cm}

\noindent {\bf Subcase (I-i)}: We assume that $1\leq j\leq b/2$ holds.\\
Note that $v_{i,j}$ and $v_{i-1,k}$ are nonadjacent in $G_{nb,b}$, where $j+1\leq k\leq b$.\ 
If $b/2+1\leq k\leq b$, then $v_{i-1,k}\notin X_{i,j}$ holds.\ 
If $j+1\leq k\leq b/2$, then we can check that $v_{i,j}\notin X_{i-1,k}$ holds.\\ 
We also note that $v_{i,j}$ and $v_{i+1,l}$ are nonadjacent, where $1\leq l\leq j-1$ and $j\ne 1$.\ 
Then $v_{i+1,l}\notin X_{i,j}$ holds.\ 


\vspace{0.2cm}
\noindent {\bf Subcase (I-ii)}: We assume that $b/2<j\leq b$ holds.\\
Note that $v_{i,j}$ and $v_{i-1,k}$ are nonadjacent in $G_{nb,b}$, where $j+1\leq k\leq b$ and $j\ne b$.\ 
Then $v_{i-1,k}\notin X_{i,j}$ holds.\\  
Note that $v_{i,j}$ and $v_{i+1,l}$ are also nonadjacent, where $1\leq l\leq j-1$.\ 
If $1\leq l\leq b/2$, then $v_{i+1,l}\notin X_{i,j}$ holds.\ 
If $b/2<l\leq j-1$, we see that $v_{i,j}\notin X_{i+1,l}$ holds.\ 


\vspace{0.2cm} 
\noindent {\bf Case (II)}: We assume that $b$ is odd.\vspace{0.2cm}\\
There is no essential difference between {\bf Case (I)} and {\bf (II)}.\ 
The difference is only the definition of the subset $X_{i,j}$ of $V(G_{nb,b})\setminus V_i$ 
by subscript.\ The reader will be able to check that 
\begin{center} 
$X_{i,1}\supseteq \cdots \supseteq X_{i,\lceil b/2\rceil -1}\supseteq X_{i,\lceil b/2\rceil }
\subseteq X_{i,\lceil b/2\rceil +1}\subseteq \cdots \subseteq X_{i,b}$
\end{center}
holds under the following definition of the subset $X_{i,j}$ of $V(G_{nb,b})\setminus V_i$ for the vertex $v_{i,j}$ of $G_{nb,b}$: 
\begin{equation*}
X_{i,j}=
\begin{cases}
N(v_{i,j})\cup \{v_{i-1,j},\ldots , v_{i-1,\lceil b/2\rceil }\}
& \text{if $1\leq j\leq \lceil b/2\rceil $,}\\
N(v_{i,j})\cup \{v_{i+1,\lceil b/2\rceil },\ldots , v_{i+1,j}\}
& \text{if $\lceil b/2\rceil <j\leq b$.}
\end{cases}
\end{equation*}
The reader also see in the same way that either $v_{i,j}\not\in X_{s,t}$ or $v_{s,t}\not\in X_{i,j}$ holds for 
any pair of nonadjacent vertices $v_{i,j}$ and $v_{s,t}$ in $G_{nb,b}$.\ 
\end{proof}
\end{thm}
We remark that a circulant graph $G_{a,b}$ has an asteroidal triple for $a\geq 3b$ and $b\geq 3$.\ 
In fact, the triple of vertices $u=1$, $v=\lceil b/2\rceil $ and $w=b$ 
becomes an asteroidal triple.\ We can find paths $P_{uv}$ with green edges, $P_{vw}$ with blue edges 
and $P_{wu}$ with red edges that we desire as in Figure 2 below.\ 
\begin{figure}[!ht]
\centering
\includegraphics[scale=1,clip]{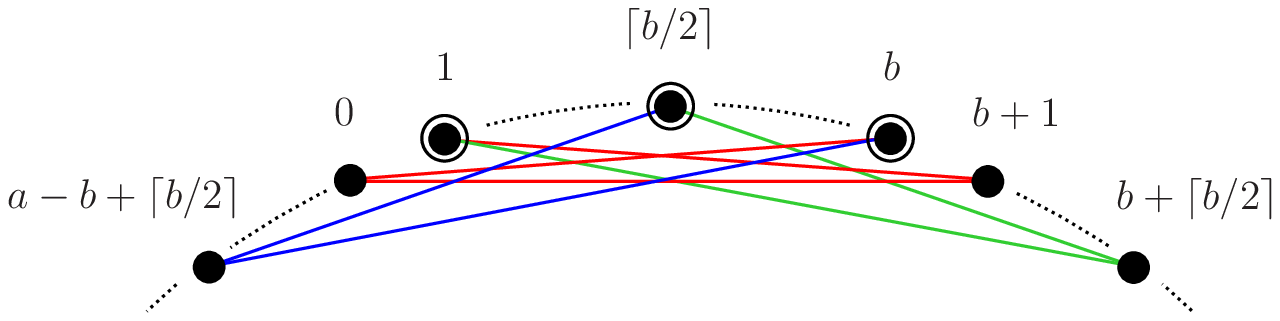}
\renewcommand{\baselinestretch}{1}
\caption[]{The graph $G_{a,b}$ with an asteroidal triple.}
\end{figure}\\
\indent Next we consider a circulant graph $G_{a,b}$ for the other cases $a=nb+r$ and $1\leq r<b$.\ 
We partition $V(G_{nb+r,b})$ into $n$ sets with consecutive $b$ vertices and a set with consecutive 
$r$ vertices so that we have an $(n+1)$-coloring of $G_{nb+r,b}$.\ 
When we take the same strategy in Theorem 4.1 for this $(n+1)$-coloring, the resulting 
family $\{X_{i,j}\}$ of sets does not satisfy condition (ii) of Theorem 3.1 except 
for the case $r=b-1$ because it has a color class with less than $b-1$ vertices.\ 
See Figure 3 and 4 below, where two vertices connected with the dashed line segment mean a pair 
of nonadjacent vertices in $G_{nb+r,b}$.\ Let $V_i$ be the 
unique color class with $r$ vertices.\ For example in Figure 3 with the assumption $r<\lceil b/2\rceil $, 
we notice that 
\begin{itemize}
\item $v_{i-1, \lceil b/2\rceil +r}$ and $v_{i+1, \lceil b/2\rceil -r}$ are nonadjacent in $G_{nb+r,b}$,  
\item $v_{i-1, \lceil b/2\rceil +r}\in X_{i+1, \lceil b/2\rceil -r}$ 
and $v_{i+1, \lceil b/2\rceil -r}\in X_{i-1, \lceil b/2\rceil +r}$ hold. 
\end{itemize}
\begin{figure}[!ht]
\centering
\includegraphics[scale=1,clip]{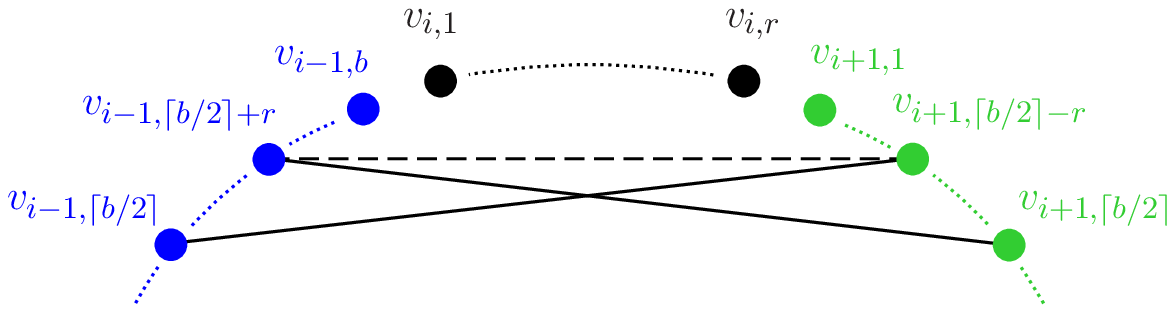}
\renewcommand{\baselinestretch}{1}
\caption[]{The graph $G_{nb+r,b}$ with $r<\lceil b/2\rceil $.}
\end{figure}
\begin{figure}[!ht]
\centering
\includegraphics[scale=1,clip]{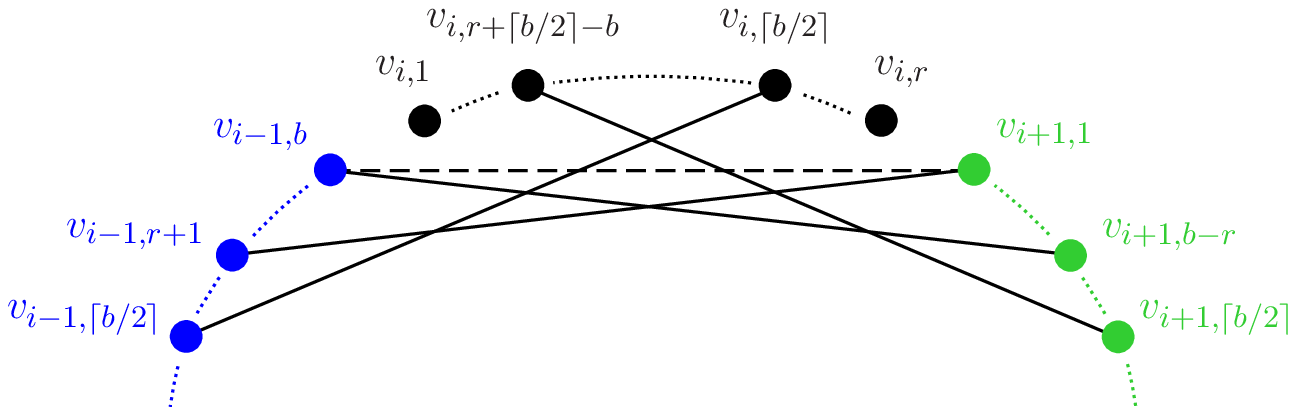}
\renewcommand{\baselinestretch}{1}
\caption[]{The graph $G_{nb+r,b}$ with $\lceil b/2\rceil \leq r<b-1$.}
\end{figure}
Also see Figure 4 with the assumption $\lceil b/2\rceil \leq r<b-1$.\ 
If $n\geq b-r-1$, we can partition $V(G_{nb+r,b})$ into $(b-r)$ sets with consecutive $(b-1)$ vertices and 
$(n-b+r+1)$ sets with consecutive $b$ vertices each of which will be a color class of an $(n+1)$-coloring of $G_{nb+r,b}$.\ 
\begin{thm}
For a circulant graph $G_{nb+r, b}$ with $n\geq 2$, $b\geq 2$ and $1\leq r<b$, 
if $n\geq b-r-1$, the inequality $\text{box}(G_{nb+r, b})\leq \chi (G_{nb+r, b})$ holds.\ 
\begin{proof}
We arrange all vertices $0,1, \ldots, nb+r-1$ of $G_{nb+r,b}$ clockwise in this order.\ 
Let $k=n-b+r+1$.\ The symbol $v_{i,j}$ means the vertex of $G_{nb+r, b}$ defined as follows: 
\begin{equation*}
v_{i,j}=
\begin{cases}
(i-1)b+j-1
& \text{if $(i,j)\in \{1,\ldots , k\}\times \{1,\ldots , b\}$,}\\
(i-1)(b-1)+j-1+k
& \text{if $(i,j)\in \{k+1, \ldots , n+1\}\times \{1,\ldots , b-1\}$.}
\end{cases}
\end{equation*}
Moreover we define 
\begin{equation*}
V_i=
\begin{cases}
\{v_{i,1}, \ldots , v_{i,b}\}
& \text{if $i\in \{1,\ldots , k\}$,}\\
\{v_{i,1}, \ldots , v_{i,b-1}\} 
& \text{if $i\in \{k+1, \ldots , n+1\}$.}
\end{cases}
\end{equation*}
We notice that $\{V_i\}_{i=1}^{n+1}$ becomes an $(n+1)$-coloring of $G_{nb+r,b}$.\  
Since $n+r/b=(nb+r)/b\leq \chi (G_{nb+r,b})$ holds, we conclude that $\chi (G_{nb+r,b})=n+1$ holds.\   

Let $|V_i|$ be the cardinality of $V_i$, that is, $|V_i|\in \{b-1, b\}$.\ 
In what follows, we always identify $V_0$ and $V_{n+2}$ with $V_{n+1}$ and $V_1$ respectively.\ 
For the vertex $v_{i,j}$ of $G_{nb+r,b}$, we define the subset $X_{i,j}$ of $V(G_{nb+r,b})\setminus V_i$ to 
be the union $N(v_{i,j})\cup Y_{i,j}$, where \vspace{0.2cm}\\
(i) when $|V_i|=b$ is even, 
\begin{equation*}
Y_{i,j}=
\begin{cases}
\{v_{i-1,j},\ldots , v_{i-1,|V_i|/2}\} & \text{if $1\leq j\leq |V_i|/2$ and $|V_{i-1}|=b$},\\
\{v_{i-1, \max \{1,j-1\}},\ldots , v_{i-1,|V_i|/2-1}\} & \text{if $1\leq j\leq |V_i|/2$ and $|V_{i-1}|=b-1$},\\
\{v_{i+1, |V_i|/2+1},\ldots , v_{i+1, \min \{j, |V_{i+1}|\}}\} & \text{if $|V_i|/2<j\leq |V_i|$},  
\end{cases}
\end{equation*}
(ii) when $|V_i|=b-1$ is even, 
\begin{equation*}
Y_{i,j}=
\begin{cases}
\{v_{i-1,j},\ldots , v_{i-1,|V_i|/2}\} & \text{if $1\leq j\leq |V_i|/2$ and $|V_{i-1}|=b$},\\
\{v_{i-1, \max \{1,j-1\}},\ldots , v_{i-1,|V_i|/2-1}\} & \text{if $1\leq j\leq |V_i|/2$ and $|V_{i-1}|=b-1$},\\
\{v_{i+1,|V_i|/2+2},\ldots , v_{i+1, \min \{j+1, |V_{i+1}|\}}\} & \text{if $|V_i|/2<j\leq |V_i|$}, 
\end{cases}
\end{equation*}
(iii) when $|V_i|=b$ is odd, 
\begin{equation*}
Y_{i,j}=
\begin{cases}
\{v_{i-1,j},\ldots , v_{i-1,\lceil |V_i|/2\rceil }\} & \text{if $1\leq j\leq \lceil |V_i|/2\rceil $ and $|V_{i-1}|=b$},\\
\{v_{i-1, \max \{1,j-1\}},\ldots , v_{i-1,\lceil |V_i|/2\rceil -1}\} & \text{if $1\leq j\leq \lceil |V_i|/2\rceil $ and $|V_{i-1}|=b-1$},\\
\{v_{i+1,\lceil |V_i|/2\rceil },\ldots , v_{i+1, \min \{j, |V_{i+1}|\}}\} & \text{if $\lceil |V_i|/2\rceil <j\leq |V_i|$,}
\end{cases}
\end{equation*}
(iv) when $|V_i|=b-1$ is odd, 
\begin{equation*}
Y_{i,j}=
\begin{cases}
\{v_{i-1,j},\ldots , v_{i-1,\lceil |V_i|/2\rceil }\} & \text{if $1\leq j\leq \lceil |V_i|/2\rceil $ and $|V_{i-1}|=b$},\\
\{v_{i-1, \max \{1,j-1\}},\ldots , v_{i-1,\lceil |V_i|/2\rceil -1}\} & \text{if $1\leq j\leq \lceil |V_i|/2\rceil $ and 
$|V_{i-1}|=b-1$},\\
\{v_{i+1,\lceil |V_i|/2\rceil +1},\ldots , v_{i+1, \min \{j+1, |V_{i+1}|\}}\} & \text{if $\lceil |V_i|/2\rceil <j\leq |V_i|$.}
\end{cases}
\end{equation*}
It is easy to check that  
\begin{center} 
$X_{i,1}\supseteq \cdots \supseteq X_{i,|V_i|/2}$ and 
$X_{i,|V_i|/2+1}\subseteq \cdots \subseteq X_{i,|V_i|}$
\end{center}
hold for a color class $V_i$ with an even number of vertices, and 
\begin{center} 
$X_{i,1}\supseteq \cdots \supseteq X_{i,\lceil |V_i|/2\rceil -1}\supseteq X_{i,\lceil |V_i|/2\rceil }
\subseteq X_{i,\lceil |V_i|/2\rceil +1}\subseteq \cdots \subseteq X_{i,|V_i|}$
\end{center}
holds for a color class $V_i$ with an odd number of vertices.\  

Note that all vertices which are nonadjacent to $v_{i,j}$ are within $V_{i-1}\cup V_i\cup V_{i+1}$.\ 
In what follows, we always assume that $v_{i,j}$ is neither adjacent to $v_{i-1,s}$ nor $v_{i+1,t}$.\ 

\vspace{0.2cm}  
\noindent {\bf Case (I)}: We assume that $j\in \{1,2,\ldots , \lceil |V_i|/2\rceil \}$. \vspace{0.2cm}\\
We may assume that $|V_{i-1}|=b$.\   
When $|V_{i-1}|=b-1$, we should replace the ordered triple $(j+1, \lceil |V_i|/2\rceil , \lceil |V_i|/2\rceil +1)$ 
in the following argument with $(j, \lceil |V_i|/2\rceil -1, \lceil |V_i|/2\rceil )$.\ 

We notice that $s\in \{j+1, j+2, \ldots , |V_{i-1}|\}$.\ 
If $j+1\leq s\leq \lceil |V_i|/2\rceil $, we see that $v_{i,j}\not\in X_{i-1,s}$ holds 
since $\lceil |V_i|/2\rceil \leq \lceil |V_{i-1}|/2\rceil  $.\ If $\lceil |V_i|/2\rceil +1\leq s\leq |V_{i-1}|$, 
then $v_{i-1,s}\not\in X_{i,j}$ holds.\ Also note that $v_{i+1,t}\not\in X_{i,j}$ holds.\ 

\vspace{0.2cm}  
\noindent {\bf Case (II)}: We assume that $j\in \{\lceil |V_i|/2\rceil +1, \ldots , |V_i| \}$ and $|V_i|$ is odd.\vspace{0.2cm}\\ 
Again first we assume that $|V_i|=b$.\ When $|V_i|=b-1$, we should replace the ordered 
triple $(j-1, \lceil |V_i|/2\rceil -1, \lceil |V_i|/2\rceil )$ in the following argument with  
$(j, \lceil |V_i|/2\rceil , \lceil |V_i|/2\rceil +1)$.\ 

We see that $t\in \{1,2,\ldots , j-1\}$.\ 
If $1\leq t\leq \lceil |V_i|/2\rceil -1$, then $v_{i+1,t}\not\in X_{i,j}$ holds.\ If $\lceil |V_i|/2\rceil \leq t\leq j-1$, 
we can check that $v_{i,j}\not\in X_{i+1,t}$ holds since $\lceil |V_{i+1}|/2\rceil \leq \lceil |V_i|/2\rceil $ holds.\ 
We also see that $v_{i-1,s}\not\in X_{i,j}$ holds.\ 

\vspace{0.2cm}  
\noindent {\bf Case (III)}: We assume that $j\in \{\lceil |V_i|/2\rceil +1, \ldots , |V_i| \}$ and $|V_i|$ is even.\vspace{0.2cm}\\ 
If $|V_i|=b$, we should just replace the ordered pair $(\lceil |V_i|/2\rceil -1, \lceil |V_i|/2\rceil )$ in Case (II) under 
$|V_i|=b$ with $(\lceil |V_i|/2\rceil , \lceil |V_i|/2\rceil +1)$.\ If $|V_i|=b-1$, we should replace the ordered triple 
$(j-1, \lceil |V_i|/2\rceil -1, \lceil |V_i|/2\rceil )$ in Case (II) under $|V_i|=b$ with 
$(j, \lceil |V_i|/2\rceil +1, \lceil |V_i|/2\rceil +2)$.\

\vspace{0.2cm}
\noindent Thus our claim follows from Theorem 3.1.\  
\end{proof}
\end{thm}
\begin{question}
Does the inequality $\text{box}(G_{a,b})\leq \chi (G_{a,b})$ hold for a circulant graph $G_{a,b}$ with $a\geq 2b$ in general?
\end{question}
For $b\geq 5$, the circulant graph $G_{2b+1,b}$ does not satisfy the condition $n\geq b-r-1$ in Theorem 4.2, but 
$\text{box}(G_{2b+1, b})\leq \chi (G_{2b+1, b})$ holds.\ In fact, $G_{2b+1,b}$ is isomorphic to a cycle 
with $2b+1$ vertices.\

 
\end{document}